\documentclass[a4paper,12pt]{amsart}

\usepackage{amssymb,pstricks,amscd}
\begin{document}
\newtheorem{cor}{Corollary}%[section]
\newtheorem{theorem}[cor]{Theorem}
\newtheorem{prop}[cor]{Proposition}
\newtheorem{lemma}[cor]{Lemma}
\theoremstyle{definition}
\newtheorem{defi}[cor]{Definition}
\theoremstyle{remark}
\newtheorem{remark}[cor]{Remark}
\newtheorem{example}[cor]{Example}

\newcommand{\aps}{\mathrm{APS}}
\newcommand{\cC}{\mathcal{C}}
\newcommand{\coker}{\mathrm{coker}}
\newcommand{\codim}{\mathrm{codim}}
\newcommand{\cone}{{\mathrm{cone}}}
\newcommand{\cun}{\cC^{\infty}}
\newcommand{\cyl}{{\mathrm{cyl}}}
\newcommand{\cz}{{\mathbb C}}
\newcommand{\Dco}{{D_{\mathrm{cone}}}}
\newcommand{\Dcy}{{D_{\mathrm{cyl}}}}
\newcommand{\Dpco}{{D^+_{\mathrm{cone}}}}
\newcommand{\Dpcy}{{D^+_{\mathrm{cyl}}}}
\newcommand{\Dmcy}{{D^-_{\mathrm{cyl}}}}
\newcommand{\Dmpco}{{D^-_{p,\mathrm{cone}}}}
\newcommand{\Dmpcy}{{D^-_{p,\mathrm{cyl}}}}
\newcommand{\Dpmco}{{D^\pm_{\mathrm{cone}}}}
\newcommand{\Dpmcy}{{D^\pm_{\mathrm{cyl}}}}
\newcommand{\Dpcos}{{D^+_{\mathrm{cone}}}^*}
\newcommand{\Dppcos}{{D^{+,*}_{p,\mathrm{cone}}}}
\newcommand{\Dppco}{{D^+_{p,\mathrm{cone}}}}
\newcommand{\Dppwco}{{D^+_{p,w,\mathrm{cone}}}}
\newcommand{\Dppcy}{{D^+_{p,\mathrm{cyl}}}}
\newcommand{\dom}{\mathrm{dom}}
\newcommand{\fm}{f^{-1}}
\newcommand{\gco}{{g_{\mathrm{cone}}}}
\newcommand{\gcy}{{g_{\mathrm{cyl}}}}
\newcommand{\ind}{\mathrm{index}}
\newcommand{\nz}{{\mathbb N}}
\newcommand{\oDp}{\overline{D^+}}
\newcommand{\oDpco}{\overline{D^+_{\mathrm{cone}}}}
\newcommand{\oX}{\overline{X}}
\newcommand{\pt}{\partial_t}
\newcommand{\px}{\partial_x}
\newcommand{\Ran}{\mathrm{Ran}}
\newcommand{\rz}{{\mathbb R}}
\newcommand{\Spec}{\operatorname{Spec}}
\newcommand{\supp}{\mathrm{supp}}
\newcommand{\tf}{{\tilde{f}}}
\newcommand{\tg}{{\tilde{g}}}
\newcommand{\tih}{\tilde{h}}
\newcommand{\td}{{\frac{t}{2}}}
\newcommand{\Tr}{\operatorname{Tr}}
\newcommand{\ud}{{\frac12}}
\newcommand{\Vol}{\operatorname{Vol}}

\title[The $L^p$ Dirac index on conical manifolds]
{On the $L^p$ index of spin Dirac operators on conical manifolds}
\author{Andr\'e Legrand}
\address{Universit\'e Paul Sabatier\\UFR MIG\\
118 route de Narbonne\\31062 Toulouse, France}
\email{legrand@picard.ups-tlse.fr}
\author{Sergiu Moroianu}
\thanks{We acknowledge the support of the
Research and Training Networks  HPRN-CT-1999-00118 ``Geometric Analysis''
and HPRN-CT-2002-00280 
``Quantum Spaces -- Noncommutative Geometry''
funded by the European Commission}
\address{Institutul de Matematic\u{a} al Academiei Rom\^{a}ne\\
P.O. Box 1-764\\RO-70700
Bucharest, Romania}
\email{moroianu@alum.mit.edu}
\date{\today}
\begin{abstract}
We compute the index of the Dirac operator on spin Riemannian
manifolds with conical singularities, acting from $L^p(\Sigma^+)$ 
to $L^q(\Sigma^-)$ with $p,q>1$. 
When $1+\frac{n}{p}-\frac{n}{q}> 0$  we obtain the usual 
Atiyah-Patodi-Singer formula, but with a spectral cut at 
$\frac{n+1}{2}-\frac{n}{q}$ instead of $0$
in the definition of the eta invariant. 
In particular we reprove Chou's formula for the $L^2$ index.
For $1+\frac{n}{p}-\frac{n}{q}\leq 0$ 
the index formula contains an extra term related to the Calder\'on projector.
\end{abstract}
\maketitle

\section{Introduction}
Let $(X,\gco)$ be a spin Riemannian manifold
 with isolated conical singularities.
The $L^2$ index of the Dirac operator $\Dpco$ on 
$X$ was computed in \cite{chou}. In this paper we
first derive Chou's formula from the index formula of 
Atiyah, Patodi and Singer \cite{aps1}. The viewpoint of this paper 
is that the passage between the two 
problems is actually elementary. We then use our method to compute the 
index of the Dirac operator from $L^p$ to $L^q$. 

\begin{theorem}\label{thmmain}
Let $(X,\gco)$ be a spin conical manifold of even dimension $n$,
and $p,q>1$. Set $\alpha_1:=\frac{n}{p}-\frac{n-1}{2}$, 
$\alpha_2:=\frac{n+1}{2}-\frac{n}{q}$.
Let
\[D^+_{p,q,\cone}:\cun_c(X,\Sigma^+)\subset L^p(X,\Sigma^+,\gco)\to 
L^q(X,\Sigma^-,\gco)\]
be the chiral Dirac operator acting as an unbounded operator
in Banach spaces.
Then for $\alpha_1+\alpha_2=1+\frac{n}{p}-\frac{n}{q}> 0$,
\begin{align*}
\ind(D^+_{p,q,\cone})&=\int_X \hat{A}(\gco)+\ud\eta_{[\alpha_2,\infty)}(M),\\
\intertext{while for $1+\frac{n}{p}-\frac{n}{q}\leq 0$,}
\ind(D^+_{p,q,\cone})&=\int_X \hat{A}(\gco)+\ud\eta_{[-\alpha_1,\infty)}(M)\\
&\quad-\dim\left(\cC_{[\alpha_2,-\alpha_1]}\right).
\end{align*}
\end{theorem}
The necessary definitions are recalled in Section \ref{secback}.
The term $\cC_{[\alpha_2,-\alpha_1]}$ is a non-negative integer
related to the Calder\'on projector. We immediately deduce Chou's result
by setting $p=q=2$. Here are other corollaries:
\begin{cor}\label{cor2}
Let $p>1$ and $\alpha:=\frac{1}{2}+\frac{n}{p}-\frac{n}{2}$, so
$\alpha> 0\Leftrightarrow p<\frac{2n}{n-1}$.
Let $p'\in(1,\infty)$ be defined by $\frac1p+\frac{1}{p'}=1$, so that 
the Dirac operator 
\[D_{p,p',\cone}:\cun_c(X,\Sigma)\subset L^p(X,\Sigma,\gco)\to L^{p'}(X,\Sigma,\gco)\]
is symmetric. 
Then for $p<\frac{2n}{n-1}$,
\begin{align*}
\ind(D^+_{p,p',\cone})&=\int_X \hat{A}(\gco)+\ud\eta_{[\alpha,\infty)}(M),\\
\intertext{while for $p\geq\frac{2n}{n-1}$,}
\ind(D^+_{p,p',\cone})&=\int_X \hat{A}(\gco)+\ud\eta_{[-\alpha,\infty)}(M)\\
&\quad-\dim\left(\cC_{[\alpha,-\alpha]}\right).
\end{align*}
\end{cor}

Next we look at the Dirac operator acting on a fixed $L^p$. Cone operators
in this context were considered by Schrohe and Seiler \cite{scse}.

\begin{cor}\label{cort}
Let $p>1$ and $\alpha:=\frac{n+1}{2}-\frac{n}{p}$. Consider the
Dirac operator
\[D^+_{p,p,\cone}:\cun_c(X,\Sigma^+)\subset L^p(X,\Sigma^+,\gco)\to 
L^p(X,\Sigma^-,\gco).\]
Then 
\[\ind(D^+_{p,p,\cone})=\int_X \hat{A}(\gco)+\ud\eta_{[\alpha,\infty)}(M).\]
\end{cor}
In particular, for $p=q=\frac{2n}{n+1}$ the index of the Dirac operator 
has the same form as the APS index formula (Theorem \ref{thaps})
 on manifolds with boundary.

Let us outline the proof of Theorem \ref{thmmain}. For simplicity consider the $L^2$ case.
First, we conjugate
$\Dco$ acting in $L^2(X,\Sigma,d\gco)$ to an unbounded operator
in $L^2(X,\Sigma,d\gcy)$, where $\gcy$ is a cylindrical
metric on $X$ conformal to $\gco$, thus transforming the problem to another 
$L^2$ index problem. Secondly, we relate this problem to the 
APS problem by restricting to a finite-length part of the cylinder.
It turns out that the $L^2$ kernel and cokernel of the conical Dirac operator
will be isomorphic to the kernel, respectively the cokernel
of an APS-type problem with a slightly different boundary spectral projection.
This is easily related to the usual APS index formula, by counting the number 
of eigenvalues between $0$ and $1/2$ of the boundary operator. Finally, 
we remark that the $\hat{A}$ form is conformally invariant.

Earlier papers on $L^2$ cohomology and index for conical manifolds
may have made the whole subject look forbiddingly technical.
We would like to reassure the reader that this is not the case here.
Except for the APS theorem which we take for granted, we give
elementary proofs of all results. In particular we do not: construct
heat kernels,
parametrices, use pseudodifferential calculi or $L^p$ Sobolev spaces 
(though the $L^2$ Sobolev embeddings on a closed manifold must be used
at one point). 

\subsection*{Acknowledgments} The second author is grateful to Andrei 
Moroianu for useful discussions. He would also like to thank the 
\emph{Equipe de G\'eom\'etrie Noncommutative de
Toulouse} for their warm hospitality at the Paul Sabatier University.

\section{Background}\label{secback}

\subsection{Conical manifolds}
The fact that $(X,\gco)$ is conical means that 
outside a compact set, $(X,\gco)$ is isometric to 
$((0,\epsilon)\times M,dr^2+r^2g_M)$, where $(M,g_M)$ is a compact, possibly 
disconnected Riemannian manifold. 

\begin{example}
Let $(\oX,g)$ be a closed Riemannian manifold and $O\in\oX$ an Euclidean point, 
in the sense that $g$ is flat in a neighborhood of $O$. In polar coordinates
we see that $\oX\setminus\{O\}$ is a conical manifold.
\end{example}
Such conical points are called \emph{fictitious}, following Dines and Schulze 
\cite{disch}. They turn out to be interesting, and we treat them 
in Section \ref{secfic}.

\subsection{Index in Banach spaces}
Let $D^+:\dom(D^+)\subset H^+\to H^-$ be a densely-defined 
unbounded linear operator between Hilbert spaces, with densely-defined adjoint
${D^+}^*$.
Let $\overline{D^+}$ be the closure of $D^+$. 
\begin{defi}\label{defind}
Assume that $\overline{D^+}$ and ${D^+}^*$ have finite-dimensional kernels.
The index of $D^+$ is defined by
\[\ind(D^+):=\dim\ker \overline{D^+}-\dim\ker {D^+}^*.\]
\end{defi} 
This is a generalization of the Fredholm index, in that the range of 
$\overline{D^+}$ needs not be closed. 

More generally, the definition holds for $H^\pm$
Banach
spaces. For completeness and since it plays a crucial role below,
we give the definition of ${D^+}^*$ and $\overline{D^+}$. 
\begin{defi}
Let $H^\pm$ be Banach spaces. 
Let $\dom({D^+}^*)$ be the space of those $u^-\in {H^-}'$
such that the map
\[\dom(D^+)\ni \phi^+\mapsto u^-(D^+\phi^+)\]
is bounded. Since $\dom(D^+)$ is assumed to be dense, this map induces
a bounded map on $H^+$. By the definition of the dual space, there exists 
a unique $u^+\in {H^+}'$ so that $u^-(D^+\phi^+)=u^+(\phi^+), \forall
\phi^+\in\dom(D^+)$. We define ${D^+}^*u^-:=u^+$.
\end{defi}

Let $H:=H^+\oplus H^-$ so $H'={H^+}'\oplus {H^-}'$.
The graph of $D^+$ is 
\[G(D^+):=\{(\phi,D^+\phi);\phi\in\dom(D^+)\}\subset H\]
and similarly for ${D^+}^*$.
Define a bilinear pairing
\begin{align*}
H\times H'&\to\rz&(\phi^+,\phi^-)\times(u^+,u^-)&\mapsto 
u^+(\phi^+)-u^-(\phi^-).
\end{align*}
Then $G(D^+)\perp G({D^+}^*)$ under this pairing. The pairing is 
continuous so $\overline{G(D^+)}\perp G({D^+}^*)$. Now we use the 
assumption that $\dom({D^+}^*)$ is dense (in practice this assumption
is checked by constructing a formal adjoint with dense domain): 
it implies that a vector of the form $(\phi^+,0)$ cannot belong to 
$\overline{G(D^+)}$ unless $\phi^+=0$. So $\overline{G(D^+)}$ is 
the graph of an operator.

\begin{defi}
$\overline{D^+}$ is the operator with closed graph $\overline{G(D^+)}$.
\end{defi}
Clearly, $\overline{D^+}^*={D^+}^*$. We define the index of $D^+$ by
Definition \ref{defind}. Note that ${D^+}^{**}=\overline{D^+}$ if and only if
$H^\pm$ are reflexive.

Specializing to the case of interest, let $X$ be a Riemannian manifold,
$\Sigma^\pm$ hermitian vector bundles, $p,q>1$ and 
\[D^+:\cun_c(X,\Sigma^+)\subset L^p(X,\Sigma^+)\to L^q(X,\Sigma^-)\]
a differential operator with formal adjoint
\[D^-:\cun_c(X,\Sigma^-)\subset L^{q'}(X,\Sigma^-)\to L^{p'}(X,\Sigma^+),\]
where $p'$ is defined by $1/p+1/{p'}=1$ and similarly for $q'$.
The elements of $L^{q'}(X,\Sigma^-)$ act as distributions on 
$L^{q}(X,\Sigma^-)$. We define the distributional action of $D^-$ on 
$u^-\in L^{q'}(X,\Sigma^-)$ in the usual way:
\[\cun_c(X,\Sigma^+)\ni\phi^+\mapsto(D^-u^-)(\phi^+):=u^-(D^+\phi^+).\]

\begin{lemma} \label{lemaddis}
A section $u^-\in L^{q'}(X,\Sigma^-)$ belongs to 
$\dom({D^+}^*)$ if and only if the distributional derivative 
$D^-u^-$ belongs to $L^{p'}(X,\Sigma^+)$. In that case, ${D^+}^*u^-=D^-u^-$.
\end{lemma}
\begin{proof}
Immediate from the definition of $\dom({D^+}^*)$ and $D^-u^-$.
\end{proof}

\subsection{The eta invariant}
Let $M$ be a compact spin manifold. For any interval $I\subset\rz$, let 
$\Pi_{I}:\cun(M,\Sigma(M))\to \cun(M,\Sigma(M))$ be the spectral projection
associated to the Dirac operator $D_M$ and $I$; more precisely, 
if $(\phi_\lambda)$ is an eigenspinor of $D_M$ of eigenvalue $\lambda$, then
\[\Pi_I(\phi_\lambda)=\begin{cases}
\phi_\lambda&\text{if $\lambda\in I$;}\\0& \text{otherwise.}\end{cases}\]
Associated to $I$ we define the complex function 
\[\eta_I(D_M,z):=\Tr\left((2\Pi_I-1)(D_M^2)^{-\frac{z}{2}}\right).\]
Special care is needed for the eigenvalue $0$, we simply define the 
complex powers to be $1$ on the nullspace of $D_M$.
The generalized eta function $\eta_I$
is well-defined and holomorphic for $\Re(z)>n-1$ and extends 
meromorphically to $\cz$ with simple poles. When $I=\rz$ we get the zeta 
function of $D_M$; when $I=[0,\infty)$ we get the extended eta function of
\cite{aps1}. The function $\eta_I$ is always regular at $z=0$ (see \cite{aps1}). 
We denote this regular(ized) value by $\eta_I(M)$. 
\begin{lemma} \label{lemeta}
For $\alpha<\beta$ we have
\[\eta_{[\alpha,\infty)}(M)=\eta_{[\beta,\infty)}(M)+2N[\alpha,\beta)\]
where $N[\alpha,\beta)$ is the number of eigenvalues of $D_M$ 
(counted with multiplicity) in the interval $[\alpha,\beta)$.
\end{lemma}
\begin{proof}
Clearly,
\[\eta_{[\alpha,\infty)}(D_M,z)-\eta_{[\beta,\infty)}(D_M,z)=
2\sum_{\substack{\alpha\leq\lambda<\beta\\ \lambda\in\Spec(D_M)}}
\lambda^{-z}.\]
Recall that we defined $0^z=1$ for all $z$. Evaluating at $z=0$ we 
get the result.
\end{proof}

\subsection{The Atiyah-Patodi-Singer index formula}
Let $Y$ be an even-dimensional compact spin manifold with boundary $M$,
and $t:Y\to(-\infty, 0]$ (the negative of) a boundary-defining function.
Fix a product decomposition $[-\epsilon,0]\times M\hookrightarrow Y$
of $Y$ near $M$, and a Riemannian metric on $Y$ which is of product type 
near the boundary:
\[g=dt^2+g_M.\] 

Over the cylinder $(-\epsilon,0]\times M$ there exist canonical 
isomorphisms of the spinor bundles $\Sigma^\pm$ with the spinor bundle of
$M$ for the induced spin structure. With these identifications, we have
\begin{align}\label{fordcyl}
\Dpcy&=\pt+D_M& \Dmcy&=-\pt+D_M
\end{align}
where $D_M$ is the Dirac operator on $M$. Set 
\[\cun(Y,\Sigma^\pm,\Pi_{I}):=\{\phi\in \cun(Y,\Sigma^\pm);\phi_{|M}\in 
\Ran(\Pi_I)\}.\]
For $\alpha\in\rz$ consider the Dirac operators with spectral boundary conditions
\begin{align*}
D^+_{\aps,\alpha}&:\cun(Y,\Sigma^+,\Pi_{[\alpha,\infty)})\to\cun(Y,\Sigma^-)\\
D^-_{\aps,\alpha}&:\cun(Y,\Sigma^-,\Pi_{(-\infty,\alpha)})\to\cun(Y,\Sigma^+)\\
\intertext{and more generally for $I\subset \rz$}
D^\pm_{\aps,I}&:\cun(Y,\Sigma^\pm,\Pi_{I})\to\cun(Y,\Sigma^\mp).
\end{align*}
These are continuous operators between Fr\'echet spaces. Being Fredholm 
in this context means having finite-dimensional kernel and co-image, which
implies that the range is closed by the open mapping theorem.

By integration by parts, \eqref{fordcyl} implies that 
$\ker D^\pm_{\aps,\alpha}$ is orthogonal to $\Ran(D^\mp_{\aps,\alpha})$
with respect to the $L^2$ inner product on $\cun(Y,\Sigma^\pm,g)$.

\begin{theorem}[\cite{aps1}]\label{thaps}
The Dirac operators $D^\pm_{\aps,0}$ are Fredholm,
\[\ker D^\pm_{\aps,0}\dotplus\Ran(D^\mp_{\aps,0})=\cun(Y,\Sigma^\pm)\]
and its index is given by
\[\ind(D^+_{\aps,0})=\int_Y\hat{A}(g)+\ud\tilde{\eta}(M),\]
where $\hat{A}(g)$ is the Hirzebruch A-hat form, and $\tilde\eta(M):=
\eta_{[0,\infty)}(M)$
is the extended eta invariant of the operator $D_M$.
\end{theorem}
The explanation for the unusual sign in front of the eta invariant is our
non-standard choice of sign for the variable $t$.

It is easy to deduce from here a similar result for $D^\pm_{\aps,\alpha}$.
\begin{cor}\label{kdep}
For $\alpha\in\rz$, the operator $D^+_{\aps,\alpha}$ is Fredholm,
\begin{equation}\label{kprec}
\ker D^\pm_{\aps,\alpha}\dotplus\Ran(D^\mp_{\aps,\alpha})=\cun(Y,\Sigma^\pm)
\end{equation}
and 
\begin{equation}\label{indalp}
\ind(D^+_{\aps,\alpha})=\int_Y\hat{A}(g)+\ud\eta_{[\alpha,\infty)}(M).
\end{equation}
\end{cor}
\begin{proof}
For simplicity assume $\alpha>0$. 
The evaluation map at the boundary
\[\cun(Y,\Sigma^+,\Pi_{[0,\infty)})/\cun(Y,\Sigma^+,\Pi_{[\alpha,\infty)})\to
\Ran(\Pi_{[0,\alpha)})\]
is an isomorphism, and by definition $\dim(\Ran(\Pi_{[0,\alpha)}))=N[0,\alpha)$.
Thus $D^+_{\aps,\alpha}$ is the restriction of $D^+_{\aps,0}$
to a Fr\'echet subspace of codimension $N[0,\alpha)$. This implies 
(using Theorem \ref{thaps}) that $D^+_{\aps,\alpha}$ has finite 
dimensional index 
\[\ind(D^+_{\aps,\alpha})=\ind(D^+_{\aps,0})-N[0,\alpha)\]
which by Lemma \ref{lemeta} implies \eqref{indalp}. Similarly
\[\ind(D^-_{\aps,\alpha})=\ind(D^-_{\aps,0})+N[0,\alpha)\]
so in particular 
\begin{equation}\label{ecuind}
\ind(D^+_{\aps,\alpha})=-\ind(D^-_{\aps,\alpha}).
\end{equation}
We have seen 
that $\ker D^\pm_{\aps,\alpha}\perp \Ran D^\mp_{\aps,\alpha}$
so
\[\begin{split}
\ind(D^+_{\aps,\alpha})&\leq \dim\ker D^+_{\aps,\alpha}
-\dim\ker D^-_{\aps,\alpha}\\
\ind(D^-_{\aps,\alpha})&\leq \dim\ker D^-_{\aps,\alpha}
-\dim\ker D^+_{\aps,\alpha}.
\end{split}\]
From \eqref{ecuind} both inequalities must be equalities, which proves
\eqref{kprec}.
\end{proof}
In the sequel we will need \eqref{kprec} in order to identify
$\codim(\Ran(D^+_{\aps,\alpha}))$ with $\dim\ker D^-_{\aps,\alpha}$.

\subsection{The Calder\'on projector}
Let $\cC\subset\cun(M,\Sigma(M))$ be the image of the Calder\'on projector,
or equivalently the space of boundary values of smooth solutions on $Y$
to the equation $D^+\phi=0$. For all $I\subset \rz$, the unique continuation property of harmonic spinors shows that the map of restriction to 
$M$ induces an isomorphism
\[\ker(D^+_{\aps,I})\to \cC\cap\Ran(\Pi_I).\]
For $J\subset I$ infinite
intervals bounded below set 
\[\cC_{I\setminus J}:=\cC\cap\Ran(\Pi_I)/\cC\cap\Ran(\Pi_J).\]
For example, $\cC_{\{0\}}=\cC\cap\Ran(\Pi_{[0,\infty)})/\cC\cap\Ran
(\Pi_{(0,\infty)})$
is the space $h_\infty$ of limiting values of extended 
$L^2$ solutions from \cite[Corollary 3.14]{aps1}. 
Clearly
\[\dim\left(\cC_{[\beta,\gamma)}\right)\leq N[\beta,\gamma).\]

Note that for $p=\frac{2n}{n-1}$, Corollary \ref{cor2} reads
\[\ind(D^+_{p,p',\cone})=\int_X \hat{A}(\gco)+\ud\eta_{[\alpha,\infty)}(M)
-h_\infty.\]

\section{Step 1: Conformal change}

To illustrate our method we first give the proof of Corollary \ref{cort}
in the case $p=2$.
The Dirac operator $\Dpmco:\cun(X,\Sigma^\pm)\to\cun(X,\Sigma^\mp)$ restricts to a
densely-defined unbounded operator 
\begin{equation}\label{dco}
\Dpmco:\cun_c(X,\Sigma^\pm)\subset L^2(X,\Sigma^\pm,d\gco)\to 
L^2(X,\Sigma^\mp,d\gco).
\end{equation}
Note that $\oDpco$ is the operator $\overline{D}_0^+$ of Chou \cite{chou}. 

Let $r:X\to(0,\infty)$ denote a smooth extension of the distance to 
the singularity. Consider the conformal metric
\begin{equation}\label{gcogcy}
\gcy:=r^{-2} \gco.
\end{equation}
The spinor bundles for the two metrics $\gcy$ and $\gco$ on $X$ are the same. 
The two Dirac operators are related by \cite[Prop. 1.3]{hitchin}:
\begin{equation}\label{hici}
\Dcy=r^{\frac{n+1}{2}}\Dco r^{-\frac{n-1}{2}}
\end{equation}
where $n$ is the dimension of $X$. Consider the isometry of Hilbert spaces
\begin{align*}
L^2(X,\Sigma, d\gco)&\to L^2(X,\Sigma, d\gcy)&
\phi&\mapsto r^{\frac{n}{2}}\phi.
\end{align*}
Using \eqref{hici}, we see that the operator $\Dpmco$ from \eqref{dco} 
is conjugated via the above isometry to
\begin{equation}\label{dcy}
r^{-\ud}\Dpmcy r^{-\ud}:\cun_c(X,\Sigma^\pm)\subset L^2(X,\Sigma^\pm,d\gcy)\to 
L^2(X,\Sigma^\mp,d\gcy).
\end{equation}
Use the change of variables 
\begin{equation}\label{chccy}
t=-\log r+\log \epsilon-1.
\end{equation}
Then $\gcy=dt^2+g_M$ for $t>-1$ so $(X,\gcy)$ is complete with 
cylindrical ends. Eq. \eqref{fordcyl} holds over the cylinder 
$(-1,\infty)\times M$. For simplicity, we define the operators
\[D^\pm:=e^{\td}\Dpmcy e^{\td}=\frac{e}{\epsilon} r^{-\ud}\Dpmcy r^{-\ud}\]
with domains given by \eqref{dcy}. Since $\Dpco$ and $\frac{\epsilon}{e}D^+$
are conjugated via an isometry, we have trivially
\begin{equation}\label{eq3}
\ind(\Dpco)=\ind(D^+).
\end{equation}

\section{Step 2: Restriction to a finite cylinder}

Let $Y\subset X$ be the compact manifold with boundary defined by $t\leq 0$.
Let $\phi\in\ker {D^+}^*$ and 
\begin{align*}
\phi(t)&=\sum_{\lambda\in\Spec D_M} a_\lambda(t)\phi_\lambda,&
D_M \phi_\lambda&=\lambda\phi_\lambda\end{align*}
the decomposition of $\phi$ over the cylinder in an orthonormal base of 
eigenspinors of $D_M$. From Lemma \ref{lemaddis}, ${D^+}^*\phi=0$ is equivalent to
$D^-\phi=0$ in distributional sense. By elliptic regularity, $\phi$ is smooth. 
Using \eqref{fordcyl}, over the cylinder $(-1,\infty)\times M$ we get
\begin{align}
&\quad(D^-\phi)(t)=0\nonumber\\
\Leftrightarrow & \quad \sum_{\lambda\in\Spec D_M} 
e^\td (-\pt (e^\td a_\lambda(t))+\lambda e^\td a_\lambda(t))\phi_\lambda=0
\nonumber\\
\Leftrightarrow &\quad a_\lambda(t)=e^{(\lambda-\ud)t}a_\lambda(0), 
\forall \lambda\in\Spec(D_M)\label{alam}
\end{align}
Since $\phi\in L^2(X,\Sigma^-,d\gcy)$ we deduce that $a_\lambda(0)=0$
for all $\lambda\geq \ud$. In other words, $\phi_{|Y}\in
\cun(Y,\Sigma^-,\Pi_{(-\infty,\ud)})$. Moreover, $\phi_{|Y}$ is a solution
of the partial differential operator $D^-$.

Conversely, let $\phi_Y\in\cun(Y,\Sigma^-,\Pi_{(-\infty,\ud)})$ be a solution
to $D^-\phi_Y=0$. Then over the cylinder $(-1,0]\times M\subset Y$, 
the coefficients of $\phi_Y$ satisfy \eqref{alam}. The spectral condition
at $t=0$ implies that $a_\lambda=0$ for $\lambda\geq\ud$.
Therefore we can consistently define $\phi\in\cun(X,\Sigma^-)$ by
\[\phi=\begin{cases}\phi_Y &\text{over $Y$;}\\ \sum_{\ud>\lambda\in\Spec D_M}
e^{(\lambda-\ud)t}a_\lambda(0)\phi_\lambda&\text{over $(-1,\infty)\times M$.}
\end{cases}\]
It is clear that $\phi$ is a distributional solution to $D^-\phi=0$ and
moreover $\phi\in L^2(X,\Sigma^-,d\gcy)$. For this last fact it is crucial that
$\lambda<\ud$. Summarizing the above discussion, we proved
\begin{prop}\label{kedps}
The map of restriction to $Y$ induces an isomorphism 
\[r_Y:\ker {D^+}^*\to\ker\left(D^-_{\aps,\ud}:\cun(Y,\Sigma^-,\Pi_{(-\infty,\ud)})
\to\cun(Y,\Sigma^+)\right).\]
\end{prop}

Similarly, let $\phi\in \ker {D^-}^*$. Then
over $(-1,\infty)\times M$ we have 
\begin{equation}\label{ptkd}
\phi(t)=\sum_{-\ud<\lambda\in\Spec D_M} e^{-(\ud+\lambda)t}a_\lambda(0)
\phi_\lambda.\end{equation}
As above, the restriction map induces an isomorphism
\[r_Y:\ker {D^-}^*\to \ker\left(D^+_{\aps,(-\ud,\infty)}:
\cun(Y,\Sigma^+,\Pi_{(-\ud,\infty)})
\to\cun(Y,\Sigma^-)\right).\] 
We are actually interested in $\ker \oDp$, an operator 
which is smaller than ${D^-}^*$, so $\ker \oDp=\ker{D^-}^*\cap \dom(\oDp)$.

\begin{lemma}\label{lemf}
Let $\phi\in \cun(X,\Sigma^+)$ satisfy \eqref{ptkd} over the cylinder 
$(-1,\infty)\times M$. Then $\phi\in\dom(\oDp)$ if and only if
$a_\lambda(0)=0$ for all $\lambda<\ud$.
\end{lemma}
\begin{proof}
First assume that $a_\lambda(0)=0$ for all $\lambda<\ud$. Let $\chi:X\to[0,1]$
be a smooth cut-off function with the properties
\[\chi(t):=\begin{cases}
1&\text{if $t<1$;}\\0& \text{if $t>2$.}
\end{cases}\]
For $k\in\nz^*$ let $\chi_k(t):=\chi(t/k)$ and define $\psi_k:=\chi_k \phi\in
\cun_c(X,\Sigma^+)$. Clearly, $\lim_{k\to\infty}\psi_k=\phi$ in $L^2$.
Moreover, 
\[D^+\psi_k=\begin{cases}
D^+\phi & \text{over $Y$;}\\\frac{1}{k} \chi'(t/k)e^t\phi &\text{over 
$(-1,\infty)\times M$}.
\end{cases}\]
The assumption on $\lambda$ and \eqref{ptkd} show that 
$|e^t\phi(t)|_{L^2(M,\Sigma(M))}$ is bounded as a function of $t$.
Since $\int_0^\infty \frac{1}{k^2}\chi'(t/k)^2dt=\frac{C}{k}\to 0$ as $k\to\infty$,
it follows that $D^+\psi_k$ converges in $L^2$ (to $D^+\phi$, a compactly 
supported smooth distribution) so $\phi\in\dom(\oDp)$.

Let now $\phi\in\dom(\oDp)$ satisfy \eqref{ptkd}. Since $\chi\phi$ has 
compact support, it follows that $(1-\chi(t))\phi\in\dom(\oDp)$. In the 
sense of distributions,
\[D^+((1-\chi(t))\phi)=-\chi'(t)e^t\phi\]
because \eqref{ptkd} implies that $D^+\phi=0$ on the support of $1-\chi$.
Thus there exists a sequence $\cun_c(X,\Sigma^+)\ni\psi_k$ such that for $k\to\infty$,
\begin{equation}\label{sec}\begin{split}
\psi_k&\stackrel{L^2}{\to}(1-\chi(t))\sum_{\substack{\lambda\in\Spec D_M\\
\lambda>-\ud}} 
e^{-(\lambda+\ud)t}a_\lambda(0)\phi_\lambda;\\
D^+\psi_k&\stackrel{L^2}{\to} {D^-}^*(1-\chi(t))\phi=-\chi'(t) 
\sum_{\substack{\lambda\in\Spec D_M\\
\lambda>-\ud}} 
e^{-(\lambda-\ud)t}a_\lambda(0)\phi_\lambda.
\end{split}\end{equation}

The right-hand side is supported in $[1,\infty)\times M$, so the
sequence $\{(1-\chi(t+1))\psi_k\}$ also fulfills \eqref{sec}.
Thus we may assume that
$\psi_k$ is supported on $[0,\infty)\times M$. 
Let $\lambda\in(-\ud,\ud)\cap\Spec(D_M)$.
We pair the second limit in \eqref{sec} with the $L^2$ distribution
\[u_\lambda=\begin{cases}0&\text{on $X\setminus(-1,\infty)\times M$;}\\
a_\lambda(0)e^{(\lambda-\ud)t}\phi_\lambda&\text{on $(-1,\infty)\times M$.}
\end{cases} \]
We get
\[\lim_{k\to\infty}(D^+\psi_k,u_\lambda)=-|a_\lambda(0)|^2\int_0^\infty 
\chi'(t) dt=|a_\lambda(0)|^2.\]
On the other hand, by definition 
\[(D^+\psi_k,u_\lambda)=(\psi_k,{D^+}^* u_\lambda)=0\]
since $D^- u_\lambda$ (in the sense of distributions) is supported at 
$t=-1$, thus outside the support of $\psi_k$. Therefore $a_\lambda(0)=0$.
\end{proof}

Thus the restriction map to $Y$ gives an isomorphism
\[r_Y:\ker \oDp\to \ker(D^+_{\aps,\ud}:\cun(Y,\Sigma^+,\Pi_{[\ud,\infty)})
\to\cun(Y,\Sigma^-)).\] 
Together with Proposition \ref{kedps} and Corollary \ref{kdep} we proved
\begin{equation}\label{eq4}
\ind(D^+)=\ind(D^+_{\aps,\ud}).
\end{equation}

\section{Step 3: The $L^2$ index formula}\label{step3}

Corollary \ref{kdep} and Lemma \ref{lemeta} give
\[\ind(D^+_{\aps,\ud})=\int_Y \hat{A}(\gcy)+\ud\tilde{\eta}(M)-N[0,\ud).\]
We claim that
\[\int_Y \hat{A}(\gcy)=\int_X \hat{A}(\gco).\]
Indeed, the $A$-hat form is a Pontryagin form; as such, it only involves the Weyl tensor 
and therefore it is conformally invariant. So $\hat{A}(\gcy)=\hat{A}(\gco)$.
Moreover, $\hat{A}(\gcy)$ vanishes on the cylinder $(-1,\infty)\times M$
by multiplicativity, so $\hat{A}(\gco)$ also vanishes on the cone
$X\setminus Y$. In particular, $\hat{A}(\gco)$ has compact support on $X$.
Consequently
\begin{align*}
\ind(\Dpco)&=\ind(D^+)&&\text{by \eqref{eq3}}\\
&=\ind(D^+_{\aps,\ud})&&\text{by \eqref{eq4}}\\
&=\int_X \hat{A}(\gco)+\ud\tilde{\eta}(M)-N[0,\ud)&&\text{by Cor.\ \ref{kdep}}.
\end{align*}
This is Chou's formula \cite[Thm. 5.23]{chou}. 

\section{The $L^p\to L^q$ index}

For $1<p\in\rz$, let $L^p(X,\Sigma^\pm,\gco)$ be the Banach space of 
$p$-integrable spinors obtained by completing $\cun_c(X,\Sigma^\pm)$
in the $L^p$ norm:
\[\|\phi\|^p_{L^p}:=\int_X|\phi|^p d\gco.\]
Let $p'\in(1,\infty)$ be the ``dual'' of $p$, i.e., 
\[\frac1p+\frac{1}{p'}=1\]
so that integration on $X$ gives a bilinear pairing 
\[L^p(X,\Sigma^\pm,\gco)\times L^{p'}(X,\Sigma^\pm,\gco)\to\cz.\]
It is well known
that $L^p$ is reflexive, so the above pairing identifies
${L^p}(X,\Sigma^\pm,\gco)'$ with $L^{p'}(X,\Sigma\pm,\gco)$.
The main result of this paper is Theorem \ref{thmmain},
the computation of the index of
\[D^+_{p,q,\cone}:\cun_c(X,\Sigma^+)\subset L^p(X,\Sigma^+,\gco)\to L^q(X,\Sigma^-,\gco)\]
in the sense of Definition \ref{defind}, for $p,q>1$. 

\begin{proof}[Proof of Theorem \ref{thmmain}]
We follow the strategy already used above for $p=q=2$. First we conjugate
$D^+_{p,q,\cone}$ to an operator acting in the cylindrical $L^p$ spaces, where 
$\gcy$ is given by Eq.\ \eqref{gcogcy}:
\[
\begin{CD}
L^p(X,\Sigma^+,\gco) @>{D^+_{p,q,\cone}}>>
&L^q(X,\Sigma^-,\gco)\\
@V{r^\frac{n}{p}\cdot}VV&@V{r^\frac{n}{q}\cdot}VV\\
L^p(X,\Sigma^+,\gcy)@>{cD^+}>> & L^q(X,\Sigma^-,\gcy)
\end{CD}
\]
with $c=(e/\epsilon)^{\alpha_1+\alpha_2}$. The vertical arrows are isometries. 
Using Eq.\ \eqref{hici} we see that 
\[\begin{split}D^+&=\left(\frac{\epsilon}{e}\right)^{\alpha_1+\alpha_2}
r^{-\alpha_2}D^+_{p,q,\cyl} r^{-\alpha_1}\\
&=e^{\alpha_2 t}D^+_{p,q,\cyl} e^{\alpha_1 t}\end{split}\]
after the coordinate change \eqref{chccy}. 
As in the $L^2$ case, it is obvious that
\begin{equation}\label{eq12}
\ind(\Dppco)=\ind(D^+).
\end{equation}
Note that the formal adjoint
\[D^-:\cun_c(X,\Sigma^-)\subset L^{q'}(X,\Sigma^-,\gcy)\to 
L^{p'}(X,\Sigma^-,\gco)\]
is given by
\[D^-=e^{\alpha_1 t}D^-_{p,q,\cyl} e^{\alpha_2 t}.\]
Let $\phi\in\ker {D^+}^*$, so by Lemma \ref{lemaddis}, 
$\phi\in L^{q'}(X,\Sigma^-,\gcy)$ is a 
distributional solution of $D^-$. By elliptic regularity, $\phi$ is smooth. 
Using Eq.\ \eqref{fordcyl}, we see that the restriction of 
$\phi$ to the cylinder $\{t\geq 0\}$ is explicitly given by the analog of 
Eq.\ \eqref{alam}:
\begin{equation}\label{ansa}
\phi(t)=\sum_{\lambda\in\Spec D_M} 
e^{(\lambda-\alpha_2)t}a_\lambda(0)\phi_\lambda.\end{equation}
\begin{lemma}\label{lemt}
Fix $w\in(1,\infty)$. A smooth spinor $\phi$ which satisfies 
\eqref{ansa} belongs to $L^w(X,\Sigma^-,\gcy)$ if and only if 
$a_\lambda(0)=0$ for all $\lambda\geq \alpha_2$.
\end{lemma}
\begin{proof}
This is clear for $w=2$ but in general it needs a proof.
First assume that $a_\lambda(0)=0$ for all $\lambda\geq \alpha_2$. Then 
the $L^2$ Sobolev norms of $\phi(t)$ decrease exponentially with $t$. 
More precisely, let $\|\phi(t)\|_{H^k(M)}:=\|D_M^k\phi(t)\|_{L^2(M)}$.
Since $\Spec D_M$ is discrete, there exists $\epsilon>0$ such that 
$\alpha_2- \lambda >\epsilon$ for all $\lambda\in\Spec D_M\cap (-\infty,\alpha_2)$.
Then 
\[\|\phi(t)\|_{H^k(M)}< C_k e^{-\epsilon t}.\]
By the Sobolev embeddings, the $C^0$ norm of $\phi(t)$ also decreases 
exponentially, so $\phi\in L^w, \forall w\geq 1$. 

Conversely, assume that
$a_\lambda(0)\neq 0$ for some $\lambda\geq \alpha_2$. We can assume 
that $a_\lambda(0)=1$. Let $C:=\|\phi_\lambda\|_{L^\infty(M,\Sigma(M))}$. Then
\begin{align*}
C\|\phi(t)\|_{L^1(M,\Sigma(M))}&\geq \int_M (\phi(t),\phi_\lambda)dx=
e^{(\lambda-\alpha_2)t}\geq 1\\
\intertext{so for $w'$ the "dual" of $w$,}
\int_M |\phi(t)|^w dx&\geq C^{-w} \Vol(M)^{-\frac{w}{w'}}.
\end{align*}
and therefore $|\phi|^{w}$ is not integrable on the cylinder. 
\end{proof}
Recall that $Y$ is the compact manifold with boundary
obtained from $X$ after removing the cylinder $(0,\infty)\times M$.
We just proved that the restriction to $Y$ defines an isomorphism
\begin{equation}\label{kdpp}
\ker {D^+}^*\simeq \ker D^-_{\aps,\alpha_2}.\end{equation}
Similarly, $\phi\in\ker {D^-}^*$ is equivalent to $\phi\in\cun(X,\Sigma^+)$,
$D^+\phi=0$ in the sense of distributions, and the restriction of $\phi$ to 
the cylinder satisfies
\begin{equation}\label{kdms}
\phi(t)=\sum_{\substack{\lambda\in\Spec D_M\\\lambda>-\alpha_1}} 
e^{-(\lambda+\alpha_1)t}a_\lambda(0)\phi_\lambda.
\end{equation}
We deduce that restriction to $Y$ gives the isomorphism
\begin{equation}\label{ekdm}
r_Y:\ker {D^-}^*\to \ker \left(D^+_{\aps,(-\alpha_1,\infty)}\right).
\end{equation}
We need to decide which elements in $\ker {D^-}^*$ live in 
$\dom(\overline{D^+})$. The series \eqref{kdms} clearly converges in $L^2$,
but we need to add to the proof of Lemma \ref{lemf} the argument for $L^p$ 
convergence. For the sake of clarity we give again the full proof.

\begin{lemma} \label{lemtt}
Let $\phi\in\cun(X,\Sigma^+)$ be of the form \eqref{kdms}
over the cylinder. Then $\phi\in\dom(\overline{D^+})$ if and only if
the coefficients of \eqref{kdms} satisfy
$a_\lambda(0)=0$ for all $\lambda<\alpha_2$.
\end{lemma}
\begin{proof}
Assume that  $a_\lambda(0)=0$ for all $\lambda<\alpha_2$.
By the Sobolev embedding theorem we prove as in Lemma \ref{lemt} that for some
$\epsilon >0$
\begin{equation}\label{decr}
|\phi(t,x)|<\begin{cases}
Ce^{-(\epsilon+\alpha_1+\alpha_2) t}&\text{for $\alpha_1+\alpha_2>0$;}\\
Ce^{-\epsilon t}&\text{for $\alpha_1+\alpha_2\leq0$.}
\end{cases}\end{equation}
so in particular $\phi\in L^p$.
Use the functions $\chi,\chi_k$ from the proof of Lemma \ref{lemf}.
The inequality \eqref{decr} shows that 
$\chi_k\phi\stackrel{L^p}{\to}\phi$. Then $D^+(\chi_k\phi)=\chi_k D^+\phi+
\chi_k'(t)e^{(\alpha_1+\alpha_2)t}\phi$. Clearly $\chi_k D^+\phi=D^+\phi$ since 
$\chi_k$ equals $1$ on the support of $D^+\phi$. Again by \eqref{decr},
the $L^q$ norm of $e^{(\alpha_1+\alpha_2) t}\phi(t)$ is bounded in $t$ and so
by changing variables,
\[\int_x|\chi_k'(t)e^{(\alpha_1+\alpha_2) t}\phi(t)|^q dx dt\leq \frac{C}{k^{q-1}}.\]
This implies that 
$D^+(\chi_k\phi)\stackrel{L^q}{\to}D^+\phi$ as $k\to\infty$. 

For the converse, there is nothing to prove if $\alpha_2\leq -\alpha_1$; 
therefore assume $\alpha_1+\alpha_2> 0$. Let $\phi \in \dom(\overline{D^+})$.
Since $\chi(t)\phi\in\dom(D^+)$ it follows that
$(1-\chi(t))\phi\in\dom(\overline{D^+})$.
By
%%%the closed graph theorem
definition, there exists a sequence $\{\psi_k\}_{k\in\nz}$ 
of compactly supported spinors such that
\begin{equation}\label{pcp}\begin{split}
\psi_k&\stackrel{L^p}{\to}(1-\chi(t))\sum_{\substack{\lambda\in\Spec D_M\\
\lambda>-\alpha_1}} 
e^{-(\lambda+\alpha_1)t}a_\lambda(0)\phi_\lambda;\\
D^+\psi_k&\stackrel{L^q}{\to} {D^-}^*((1- \chi(t))\phi)=-\chi'(t) \sum_{\substack{\lambda\in\Spec D_M\\
\lambda>-\alpha_1}} 
e^{-(\lambda-\alpha_2)t}a_\lambda(0)\phi_\lambda.
\end{split}\end{equation}
The sequence $\{(1-\chi(t+1))\psi_k\}$ also satisfies \eqref{pcp}, so we can 
assume that $\supp(\psi_k)\subset[0,\infty)\times M$.

For $\lambda\in (-\alpha_1,\alpha_2)$ consider the distribution
\[u_\lambda=\begin{cases}0&\text{on $X\setminus(-1,\infty)\times M$;}\\
a_\lambda(0)e^{(\lambda-\alpha_2)t}\phi_\lambda&\text{on $(-1,\infty)\times M$.}
\end{cases} \]
Since $\lambda<\alpha_2$, it follows that $u_\lambda\in L^p$ so the second
limit in \eqref{pcp} commutes with the pairing with $u_\lambda$:
\[\lim_{k\to\infty}(D^+\psi_k,u_\lambda)=-|a_\lambda(0)|^2\int_0^\infty 
\chi'(t) dt=|a_\lambda(0)|^2.\]
On the other hand, by definition 
\[(D^+\psi_k,u_\lambda)=(\psi_k,{D^+}^* u_\lambda)=0\]
since $D^- u_\lambda$ (in the sense of distributions) is supported at 
$t=-1$, thus outside the support of $\psi_k$. Therefore $a_\lambda(0)=0$.
\end{proof}

By Lemma \ref{lemtt} and \eqref{ekdm}, 
\begin{equation}\label{kerp}\ker \overline{\Dppco}\simeq\begin{cases}
\ker(D^+_{\aps,\alpha_2})&\text{for $\alpha_1+\alpha_2>0$;}\\
\ker(D^+_{\aps,(-\alpha_1,\infty)})& \text{for $\alpha_1+\alpha_2\leq 0$.}
\end{cases}
\end{equation} 
In the second case,
\[\dim\ker \overline{D^+}=\dim\ker D^+_{\aps,\alpha_2}-\cC_{[\alpha_2,-\alpha_1]}\]
because by definition
$\ker {D^+_{\aps,(-\alpha_1,\infty)}}\simeq \cC_{(-\alpha_1,\infty)}$.
Recall now from Section \ref{step3} that the $\hat{A}$ form 
is a conformal invariant and vanishes near the singularity.
These facts together with \eqref{kdpp}, \eqref{eq12} and Corollary \ref{kdep} finish 
the proof of Theorem \ref{thmmain}.
\end{proof}

\section{Fictitious conical singularities}\label{secfic}

Let $(\oX,g)$ be a closed spin manifold which contains a finite set
$\{O\}$ of Euclidean points, in the sense that each $O_j\in\{O\}$ has a 
flat neighborhood. Writing $g$ in polar coordinates near $O_j$, we see that 
$X:=\oX\setminus \{O\}$ is a conical spin manifold, and the basis of the 
cone is a disjoint union of spheres with the standard metric. 
For $n\geq 3$ the sphere $S^{n-1}$ has a unique spin structure, 
while for $n=2$ the spin structure 
on each circle must be bounding (non-trivial). The eigenvalues of the 
associated Dirac operator $D_{S^{n-1}}$ are
\[\pm \left(\frac{n-1}{2}+k\right), k=0,1,2,\ldots\]
with multiplicity $2^{[\frac{n-1}{2}]}\binom{k+n-2}{k}$ 
(see e.g., \cite{baersph}). The eigenspinors 
for the smallest eigenvalues $\pm \frac{n-1}{2}$ are simply restrictions
of parallel (i.e., constant) positive, respectively negative spinors from $\rz^n$. 
In particular there is no eigenvalue between $0$ and $\ud$. 
So for the $L^2$ index, Chou's formula reproved in Section \ref{step3} 
and the Atiyah-Singer formula give
\[\ind_{L^2}(\Dpco)=\ind(D^+_{\oX}).\]
This equality may come as a surprise, since the domains of the two operators are 
not the same. Moreover, the similar equality is not true for the index of 
$D^+_{p,q,\cone}$! This can be seen by simply comparing the index formulae.
We discuss below the case of $D^+_{p,p',\cone}$.

\begin{prop} 
For $p\geq\frac{n}{n-1}$
every spinor in the nullspace of $\overline{D^+_{p,p',\cone}}$, 
respectively ${D^{+,*}_{p,p',\cone}}$ 
extends to a harmonic spinor on $\oX$. Conversely, the restriction
of every harmonic spinor on $\oX$ to $X$ belongs to 
the nullspace of $\overline{D^+_{p,p',\cone}}$, respectively $D^{+,*}_{p,p',\cone}$.
\end{prop}
\begin{proof}
First notice that $\alpha_1=\alpha_2=:\alpha$ satisfies
\[-\frac{n-1}{2}<\alpha<\frac{n+1}{2}.\]
The only eigenvalue of $D_{S^{n-1}}$ in this interval is $\frac{n-1}{2}$, with 
multiplicity $2^{[\frac{n-1}{2}]}$.
Then from \eqref{kerp}, 
\begin{align*}
\ker \overline{D^{+}_{p,p',\cone}}&\simeq\begin{cases}
\ker(D^+_{\aps,\frac{n-1}{2}}) &\text{if $\alpha\leq\frac{n-1}{2}$;}\\
\ker(D^+_{\aps,\frac{n+1}{2}}) &\text{if $\alpha>\frac{n-1}{2}$;}
\end{cases}\\
\intertext{and from \eqref{kdpp},}
\ker D^{+,*}_{p,p',\cone}&\simeq\begin{cases} 
\ker(D^-_{\aps,(-\infty,-\frac{n-1}{2}]})&\text{if $\alpha\leq\frac{n-1}{2}$;}\\
\ker(D^-_{\aps,(-\infty,\frac{n-1}{2}]})&\text{if $\alpha >\frac{n-1}{2}$.}
\end{cases}
\end{align*}
Assume now that $\alpha\leq\frac{n-1}{2}$, or equivalently $p\geq\frac{n}{n-1}$.
Let $\phi_Y$ be a harmonic spinor in $\ker(D^+_{\aps,\frac{n-1}{2}})$.
From \eqref{kdms}, 
\[\phi_Y(t)=\sum_{\substack{\lambda\in\Spec D_{S^{n-1}}\\
\lambda\geq\frac{n-1}{2}}} 
e^{-(\lambda+\alpha)t}\phi_\lambda\]
for $-1<t\leq 0$, where $\phi_\lambda$ is an eigenspinor of
eigenvalue $\lambda$ but not necessarily of $L^2$-length $1$. 
Extend $\phi_Y$ to $X$ by the same formula for
$t>0$ and then pull it back to a spinor $\phi_\cone$ on the cone via the 
isometry from the proof of Theorem \ref{thmmain}. Therefore
\begin{align*}
\phi_\cone(r)&=r^{-\frac{n}{p}}\sum_{\substack{\lambda\in\Spec D_{S^{n-1}}\\
\lambda\geq\frac{n-1}{2}}} 
e^{-(\lambda+\alpha)t}\phi_\lambda\\
&=\sum_{k=0}^\infty r^k\phi_{k+\frac{n-1}{2}}
\end{align*}
so $\phi_\cone$ extends to a $L^\infty$ spinor $\phi$ on $\oX$. 
Now $\phi_\frac{n-1}{2}$ is the restriction of a \emph{constant} positive 
spinor from $\oX$ to $S^{n-1}$, so it extends smoothly in $r=0$. 
From \eqref{hici}, $D^+_{\oX}\sum_{k=1}^\infty r^k\phi_{k+\frac{n-1}{2}}$
is also in $L^\infty(\oX,\Sigma^-)$. 
Therefore the distribution $D^+_{\oX}\phi$ is on one hand in $L^\infty(\oX)$
and on the other hand it vanishes on $X$. It follows that $\phi$ is 
a solution to $D^+_{\oX}$. 

Conversely, every solution $\phi$ to $D^+_{\oX}$ restricts to a 
$L^p$ distributional solution $\phi_\cone$ to 
$\Dpco$ on $X$. Eq. \eqref{ekdm}, Lemma \ref{lemtt} and the 
condition on $\alpha$ show that $\phi_\cone$ actually belongs to
$\dom(\overline{D^{+}_{p,p',\cone}})$.

The statement about $D^{+,*}_{p,p',\cone}$ is proved in the same way.
\end{proof}

In conclusion, for $p\geq\frac{n}{n-1}$ the $L^p$ index problem on $X$
reduces to the usual index problem on $\oX$. 
For $1<p<\frac{n}{n-1}$ the eigenvalue $\frac{n-1}{2}$ of the Dirac operator on
the sphere creeps into the picture, so the $L^p$ index of $D^+$ on $X$
is $2^{[\frac{n-1}{2}]}$ less than $\hat{A}(\oX)$.

\section{Possible extensions}
Our goal was to give the simplest possible solution to the $L^p$ index problem,
so we did not cover metrics which are only asymptotically
conical. Our elementary method clearly breaks down in this case, and more 
conceptual approaches are needed, like parametrices for cone operators. 
Such parametrices can be constructed via either Melrose's $b$-calculus 
\cite{melaps} or the cone calculus of Schulze \cite{schulze}. General
elliptic cone operators in $L^p$ spaces are treated in \cite{scse}. 
But an index formula generalizing Theorem \ref{thmmain} is still missing
as of writing of this paper.

We have not discussed Fredholmness of our operators for the same reason, 
but note that the Fredholm property of $D^{+}_{p,p,\cone}$ 
seems to follow from the results of \cite{scse}.

\bibliographystyle{amsplain}

\begin{thebibliography}{9}  

\bibitem{aps1}
M.F.~Atiyah, V.K.~Patodi and I.M.~Singer,  
{\sl Spectral asymmetry and Riemannian geometry. I,} 
Math. Proc. Cambridge Philos. Soc. {\bf 77} (1975), 43--69.

\bibitem{baersph}
C.~B\"ar,
{\sl The Dirac Operator on Space Forms of Positive Curvature, } 
J. Math. Soc. Japan {\bf 48} (1996), 69--83.

\bibitem{chou}
A.~W.~Chou, 
{\sl The Dirac operator on spaces with conical singularities and positive 
scalar curvatures, }
Trans. Amer. Math. Soc. {\bf 289} (1985), no. 1, 1--40.

\bibitem{disch}
N.~Dines and B.~W.~Schulze,
{\sl Mellin-Edge Representation of Elliptic Operators, }
preprint (2003).

\bibitem{hitchin}
N.~Hitchin,
{\sl Harmonic spinors, }
Adv. in Math. {\bf 14} (1974), 1--55.

\bibitem{melaps}
R.~B.~Melrose,
{\sl The Atiyah-Patodi-Singer index theorem, }
Research Notes in Mathematics {\bf 4}, 
A. K. Peters, Wellesley, MA (1993).

\bibitem{scse}
E.~Schrohe, and J.~Seiler,
{\sl Ellipticity and invertibility in the cone algebra on $L\sb p$-Sobolev 
spaces, } Integr. Equat. Oper. Th. {\bf 41} (2001), 93--114.

\bibitem{schulze}
B.~W.~Schulze, 
{\sl Pseudo-differential boundary value problems, 
conical singularities, and asymptotics, }
Mathematical Topics {\bf 4}, Akademie Verlag, Berlin (1994).

\end{thebibliography}

\end{document}